\documentclass{amsart}
\usepackage{t1enc}
\usepackage[utf8]{inputenc}
\usepackage[english]{babel}
\usepackage{amsmath, amssymb, amsfonts, amsthm}
\usepackage[colorlinks=true]{hyperref}
\usepackage[capitalize, noabbrev]{cleveref}
\usepackage{graphicx, color, calc}
\usepackage[msc-links]{amsrefs}

\newtheorem{theorem}{Theorem}
\newtheorem{lemma}[theorem]{Lemma}
\theoremstyle{remark}
\newtheorem{remark}{Remark}

\DeclareMathOperator{\lk}{lk}

\newcommand{\Z}{\mathbb Z}
\newcommand{\Q}{\mathbb Q}

\renewcommand{\d}{\partial}
\renewcommand{\S}{\Sigma}
\newcommand{\nc}{\not\equiv}
\newcommand{\iso}{\sim_4}
\newcommand{\niso}{\not\sim_4}
\renewcommand{\v}{\mathbf v}
\newcommand{\T}{\mathcal T}

\begin{document}

\title{Non-isotopic Seifert surfaces in the 4-ball}
\author[Zs. Fehér]{Zsombor Fehér}
\address{Mathematical Institute, University of Oxford, Andrew Wiles Building, Radcliffe Observatory Quarter, Woodstock Road, Oxford, OX2 6GG, UK}
\email{feher@maths.ox.ac.uk}
\begin{abstract}
We give families of knots and links with pairs of Seifert surfaces that are topologically non-isotopic in $D^4$. This generalizes the main example of Hayden--Kim--Miller--Park--Sundberg and the proof is similarly based on the double branched cover and the Seifert form. Moreover, using Conway's theory of topographs, we implement an algorithm that decides whether two integral quadratic forms in two variables are isomorphic over $\Z$.
\end{abstract}
\maketitle

\section*{Introduction}
In 1982, Livingston \cite{Livingston} asked whether there exists a knot with two Seifert surfaces of the same genus that are not isotopic in $D^4$. This was answered positively in 2022 by Hayden--Kim--Miller--Park--Sundberg \cite{HKMPS}. Their main example (based on Lyon's construction \cite{Lyon-example}) comes from attaching a band to two annuli bounded by the torus link $T_{4,-6}$. They showed that the symmetrized Seifert form obstructs these surfaces from being topologically isotopic, and that the induced cobordism maps in Khovanov homology obstructs smooth isotopy.

Here we study when the symmetrized Seifert form can be used to distinguish Seifert surfaces that are modifications of this example, which can be denoted by $\S_i(2,-3,0,-1)$ in our 4-parameter family.

In \cref{sec:1}, we define a family of pairs of Seifert surfaces $\S_i(p,q,k,n)$, $i=0,1$, and prove that they are topologically non-isotopic in $D^4$ under certain conditions on $(p,q,k,n)$ (\cref{thm:pqkn,thm:n=0}). We also see many instances where the surfaces in this family are topologically isotopic (\cref{thm:isotopic}).  A related, more elegant family of Seifert surfaces that are bounded by 2-component links is given in \cref{thm:pq,thm:cables}. These results extend the currently known examples of non-isotopic Seifert surfaces.

In \cref{sec:2}, we summarize Conway's theory of topographs, a diagrammatic way to view quadratic forms in two variables, and give an algorithm that decides whether two such quadratic forms are isomorphic over $\Z$. We use this algorithm to find further examples of non-isotopic surface-pairs $\S_i(p,q,k,n)$ that are not covered by our previous theorems. Furthermore, the structure of the topographs also explains the various conditions in our theorems.

\section{Families of surface-pairs}\label{sec:1}

\subsection{Main results}

We will use the notation $F_0\iso F_1$ for surfaces $F_i\subset S^3$ that are topologically isotopic relative to boundary in $D^4$.

\begin{theorem}\label{thm:pq}
Let $p,q$ be coprime integers. The torus link $T_{2p,2q}$ cuts the torus into two annuli, giving two Seifert surfaces $A_0(p,q)$ and $A_1(p,q)$ for $T_{2p,2q}$. Then $A_0(p,q)\niso A_1(p,q)$ if and only if $|p|,|q|>1$.
\end{theorem}

Note that if $p,q$ are not coprime, then we can still cut the torus into two surfaces bounded by $T_{2p,2q}$, but these surfaces are not connected, and hence are trivially non-isotopic as they connect different component-pairs of $T_{2p,2q}$.

This theorem will be proven as a consequence of \cref{thm:pqkn}. There we investigate whether the surfaces $A_0(p,q)$ and $A_1(p,q)$ remain non-isotopic if we attach a common band $B\subset S^3$ to their boundary, giving two Seifert surfaces for a knot. Note that if $F_0\iso F_1$, then $F_0\cup B\iso F_1\cup B$, because the isotopy between the $F_i$ can be done `above' the band $B$ in the fourth dimension. However, the opposite is not true: $F_0\niso F_1$ does not imply $F_0\cup B\niso F_1\cup B$. We will show that for a family of bands $B(k,n)$, the surfaces $A_i(p,q)\cup B(k,n)$ for $i=0,1$ are non-isotopic. Then it follows that $A_i(p,q)$ are non-isotopic as well.

\begin{figure}[htb]
\par{\centering\def\svgwidth{0.95\textwidth}
\begingroup%
  \makeatletter%
  \providecommand\color[2][]{%
    \errmessage{(Inkscape) Color is used for the text in Inkscape, but the package 'color.sty' is not loaded}%
    \renewcommand\color[2][]{}%
  }%
  \providecommand\transparent[1]{%
    \errmessage{(Inkscape) Transparency is used (non-zero) for the text in Inkscape, but the package 'transparent.sty' is not loaded}%
    \renewcommand\transparent[1]{}%
  }%
  \providecommand\rotatebox[2]{#2}%
  \newcommand*\fsize{\dimexpr\f@size pt\relax}%
  \newcommand*\lineheight[1]{\fontsize{\fsize}{#1\fsize}\selectfont}%
  \ifx\svgwidth\undefined%
    \setlength{\unitlength}{816.17051012bp}%
    \ifx\svgscale\undefined%
      \relax%
    \else%
      \setlength{\unitlength}{\unitlength * \real{\svgscale}}%
    \fi%
  \else%
    \setlength{\unitlength}{\svgwidth}%
  \fi%
  \global\let\svgwidth\undefined%
  \global\let\svgscale\undefined%
  \makeatother%
  \begin{picture}(1,0.41708204)%
    \lineheight{1}%
    \setlength\tabcolsep{0pt}%
    \put(0,0){\includegraphics[width=\unitlength,page=1]{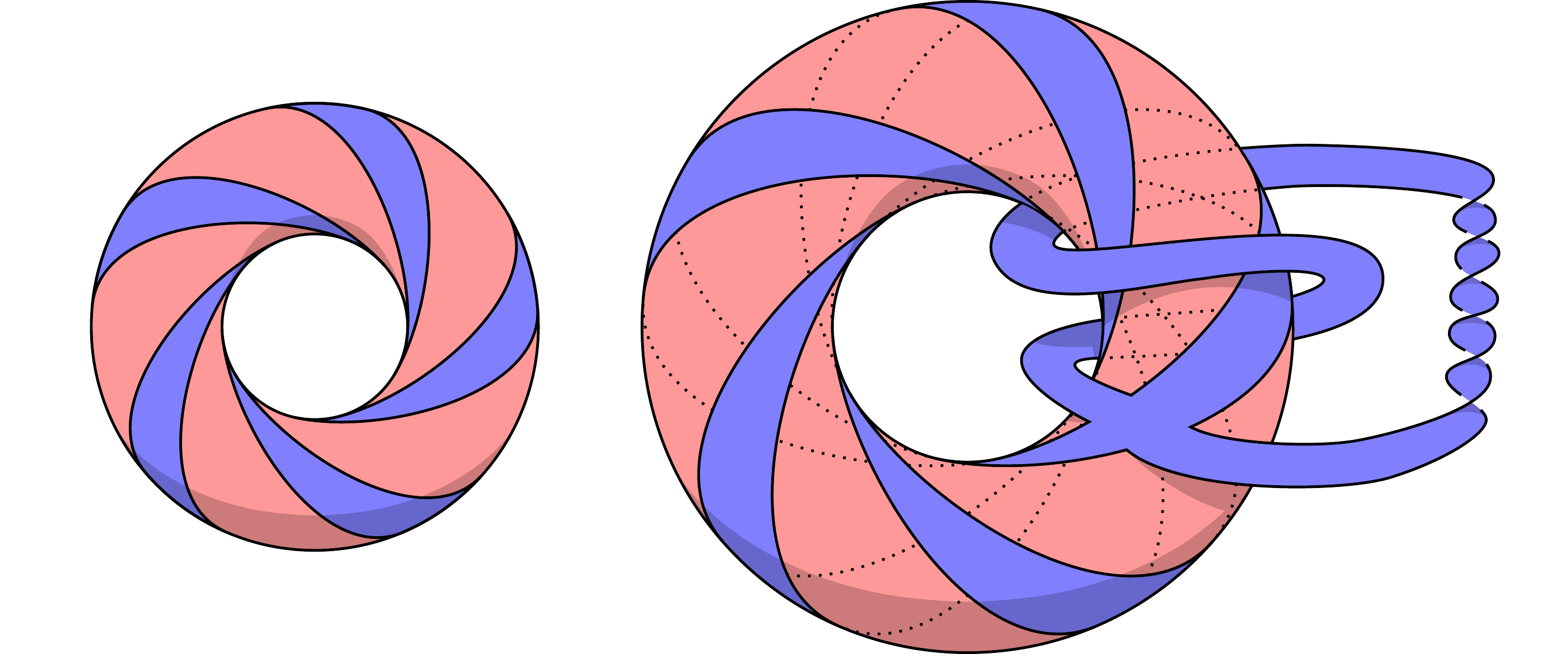}}%
    \put(0.00188473,0.33387761){\color[rgb]{0.89019608,0,0}\makebox(0,0)[lt]{\lineheight{1.25}\smash{\begin{tabular}[t]{l}$A_0(p,q)$\end{tabular}}}}%
    \put(-0.00111276,0.07896663){\color[rgb]{0,0,0.89019608}\makebox(0,0)[lt]{\lineheight{1.25}\smash{\begin{tabular}[t]{l}$A_1(p,q)$\end{tabular}}}}%
    \put(0.79837534,0.36657467){\color[rgb]{0,0,0.89019608}\makebox(0,0)[lt]{\lineheight{1.25}\smash{\begin{tabular}[t]{l}$\S_1(p,q,k,n)$\end{tabular}}}}%
    \put(0.62409587,0.20273636){\makebox(0,0)[lt]{\lineheight{1.25}\smash{\begin{tabular}[t]{l}$k$\end{tabular}}}}%
    \put(0.9730209,0.22070087){\makebox(0,0)[lt]{\lineheight{1.25}\smash{\begin{tabular}[t]{l}$n$\end{tabular}}}}%
    \put(0.82318905,0.0704291){\makebox(0,0)[lt]{\lineheight{1.25}\smash{\begin{tabular}[t]{l}$B(k,n)$\end{tabular}}}}%
  \end{picture}%
\endgroup%
\par}
\caption{The surfaces in the case $(p,q)=(3,5)$, $k=2$, and $n=3$.}\label{fig:pqkn1}
\end{figure}
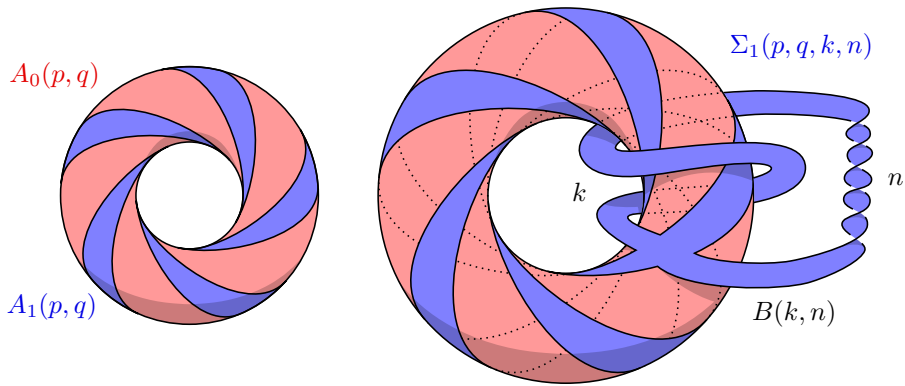

Let $(p,q,k,n)$ be integers with $p,q$ coprime. Consider a band $B(k,n)$ as shown in \cref{fig:pqkn1}, that (viewed as a framed oriented closed curve) has linking number $k$ with the core of the solid torus and has self-linking number~$n$. Let $\S_i(p,q,k,n)=A_i(p,q)\cup B(k,n)$. Note that $B(k,n)$ can refer to any knotted band with linking numbers $k$ and $n$, so $\S_i(p,q,k,n)$ ($i=0,1$) actually refers to any member of an infinite family of pairs of surfaces.

\begin{theorem}\label{thm:pqkn}
Let $p,q>1$ be coprime, $k$ an integer such that $2kp\nc1 \pmod q$, and $n$ an integer with
\[
n \ge \frac{k(pk-1)}q + \left( \frac{pq}{12} - \frac16 + \frac{1}{2pq} \right).
\]
Then after pushing the interior of $\S_i(p,q,k,n)=\S_i$ into $D^4$, the pairs $(D^4,\S_0)$ and $(D^4,\S_1)$ are not homeomorphic. In particular, $\S_0\niso\S_1$.
\end{theorem}

We will see that in some cases the lower bound on $n$ can be improved, e.g.~if $k=0$, or if $q \le 3$, then $n > \frac{k(pk-1)}q$ is sufficient. However, the next theorem shows that we cannot go lower than that.

\begin{theorem}\label{thm:isotopic}
Let $p,q>1$ be coprime, and suppose that integers $k,n$ satisfy
\[
n=\frac{k(pk-1)}q.
\]
Then $\S_0(p,q,k,n)\iso\S_1(p,q,k,n)$.
\end{theorem}

\begin{theorem}\label{thm:n=0}
Let $p,q>1$ be coprime and $n=0$. Then there exists an integer $k_0$ such that if $|k|>k_0$ is an integer satisfying $2kp\nc 1 \pmod q$, then after pushing the interior of $\S_i(p,q,k,n)=\S_i$ into $D^4$, the pairs $(D^4,\S_0)$ and $(D^4,\S_1)$ are not homeomorphic. In particular, $\S_0\niso\S_1$.
\end{theorem}

See \cref{sec:discussion} for more discussion about the conditions in \cref{thm:pqkn,thm:isotopic,thm:n=0}. In particular, the diagrams in \cref{fig:parabolas} display these conditions very clearly.

We can generalize \cref{thm:pq} for 2-component cables in the following way.

\begin{theorem}\label{thm:cables}
Let $p,q>1$ be coprime and $K$ be any knot. Let $L$ be the $(2p,2q)$-cable of $K$. Then $L$ has non-isotopic Seifert surfaces. Namely, $L$ cuts the boundary of a tubular neighbourhood of $K$ into two annuli, giving two Seifert surfaces $A_0$ and $A_1$ for $L$. Then $A_0\niso A_1$.
\end{theorem}

In fact, just like we allowed $B(k,n)$ to be any knotted band in the definition of the surface-pairs $\S_i(p,q,k,n)$, we can also allow $A_0(p,q)\cup A_1(p,q)$ to be any knotted torus (bounding a solid torus disjoint from the inside of $B(k,n)$), and \cref{thm:pqkn,thm:isotopic,thm:n=0} remain true. We will see this in the proof of \cref{thm:cables}.

\subsection{Proofs}

To prove \cref{thm:pqkn}, we will need the following lemma.

\begin{lemma}\label{thm:uv0v1}
Let $u,v_0,v_1$ be integers such that $v_0\nc \pm v_1 \pmod u$. For an integer~$t$, consider the quadratic forms
\[
Q_i(x,y)=(ux+v_iy)^2+ty^2.
\]
Then for all
\[
t > t_0 = \frac{(u/2-1)^2}{3},
\]
the quadratic forms $Q_0, Q_1$ are not isomorphic over $\Z$.
\end{lemma}

\begin{remark}\label{rem:t0}
The lower bound on $t$ might be weakened, but at the cost of having a more cumbersome, number-theoretic formula. For an integer $v$, let us denote by $[v]_u$ the integer that satisfies $0\le[v]_u\le u/2$ and $[v]_u\equiv \pm v\pmod u$. If $[v_0]_u<[v_1]_u$ (swap $v_0$ and $v_1$ if not), then this formula is
\[
t > t_1 = \max_{c\ge 2\text{ integer}} \frac{[v_0]_u^2-[cv_1]_u^2}{c^2-1}.
\]
\end{remark}

\begin{remark}\label{rem:quad11}
The condition on $t$ cannot be omitted: e.g.~the quadratic forms
\[
(30x+7y)^2+ty^2 \quad \text{and} \quad (30x+13y)^2+ty^2
\]
are not isomorphic for $t\ge 12$, but are isomorphic for $t=11$. The formula in the previous Remark gives $t_1=11$, while $t_0=65+\frac13$.
\end{remark}

\begin{proof}
We can assume $|v_i|\le u/2$, as taking $Q_i(x\pm y, y)$ (which is isomorphic to $Q_i$) increases/decreases the value of $v_i$ by $u$. We can also assume $0\le v_i$, as taking $Q_i(x, -y)$ replaces $v_i$ by $-v_i$. These manipulations also do not change the condition $v_0\nc \pm v_1 \pmod u$. Hence, without loss of generality, we can assume $0\le v_0<v_1\le u/2$.

Let $H_i$ denote the set of integers $(x,y)$ for which
\[Q_i(x,y)=v_0^2+t.\]
If $t>0$, then $Q_i$ is positive definite, hence $H_i$ is finite. We will show $|H_0|>|H_1|$, which proves that $Q_0$ and $Q_1$ are not isomorphic.

For $y=0$, the values $(x,y)$ contribute the same to $H_0$ and $H_1$, as $Q_0(x,0)=Q_1(x,0)$.
The case $|y|=1$ gives the elements $(0,\pm1)\in H_0$. However, it gives no elements in $H_1$, as $v_0<v_1\le u/2$ implies $Q_1(x,\pm1)>v_0^2+t$ for all $x$.

We show that for $t$ large enough, there are no other elements in $H_1$. As above, $[v]_u$ denotes the integer that satisfies $0\le[v]_u\le u/2$ and $[v]_u\equiv \pm v\pmod u$. If $|y|=c\ge 2$, then the smallest value that $Q_1(x,y)$ can take is $[cv_1]_u^2+tc^2$. Hence it suffices to have
\begin{align*}
[cv_1]_u^2+tc^2 &> v_0^2+t, \\
t &> \frac{v_0^2-[cv_1]_u^2}{c^2-1}
\end{align*}
for all integers $c\ge 2$. Since $v_0\le v_1-1\le u/2-1$ and $[cv_1]_u^2\ge 0$,
\[
\frac{v_0^2-[cv_1]_u^2}{c^2-1}\le \frac{(u/2-1)^2}{2^2-1}=t_0
\]
is suitable. If $t>t_0$, then $|H_0|>|H_1|$, and so $Q_0$ and $Q_1$ are not isomorphic over $\Z$.
\end{proof}

\begin{figure}[ht]
\par{\centering\def\svgwidth{1\textwidth}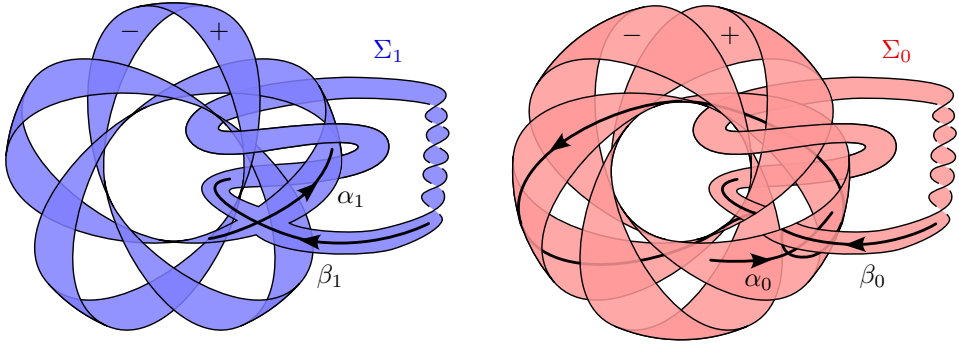\par}
\caption{Orientations for the base of $H_1(\S_i)$, $(r,s)=(1,2)$ pictured.}\label{fig:pqkn2}
\end{figure}

\begin{proof}[Proof of \cref{thm:pqkn}]
Let $X_i$ be the double branched cover of $D^4$ along $\S_i$. Then the intersection form of $H_2(X_i)$ is given by the symmetrized Seifert form of $\S_i\subset S^3$ \cite{Gordon}.

A base for $H_1(\S_1)$ is given by $\alpha_1,\beta_1$ (see \cref{fig:pqkn2}), where $\alpha_1$ is the core of $A_1(p,q)$ and $\beta_1$ is the core of the band $B(k,n)$ with opposite orientation and its endpoints joined by a small arc in $A_1(p,q)$. Then $\alpha_1$ is a torus knot $T_{p,q}$ and has self-linking number $pq$. Thus, the Seifert matrix of $\S_1$ is
\[
V_1=
\begin{bmatrix}
pq   & -kp\\
1-kp & n
\end{bmatrix}.
\]

A base for $H_1(\S_0)$ is $\alpha_0,\beta_0$, where $\alpha_0$ is a torus knot $T_{p,q}$ in $A_0(p,q)$, but for $\beta_0$ we join together the endpoints of the band $B(k,n)$ by some longer arc in $A_0(p,q)$. Then $\beta_0-\beta_1$ (meant as a union of oriented curves, cancelling the oppositely oriented overlapping parts) is a closed curve on the torus, hence it is a torus knot $T_{r,s}$. Since $T_{r,s}$ intersects $T_{p,q}$ once positively, the integers $r,s$ satisfy
\[
ps-qr=1.
\]
Using $\beta_0=T_{r,s}+\beta_1$, we get
\begin{align*}
\lk(\beta_0, \alpha_0^+)
&=\lk(T_{r,s}, \alpha_0^+) + \lk(\beta_1, \alpha_0)
=ps-kp, \\
\lk(\beta_0,\beta_0^+)
&=\lk(T_{r,s},T_{r,s}^+)+2\lk(T_{r,s},\beta_1^+)+\lk(\beta_1,\beta_1^+)
=rs-2kr+n.
\end{align*}
Therefore, the Seifert matrix of $\S_0$ is
\[
V_0=
\begin{bmatrix}
pq    & qr-kp\\
ps-kp & rs-2kr+n
\end{bmatrix}.
\]

The quadratic forms associated to the symmetrized Seifert matrices are
\begin{align*}
Q_1(x,y) &=
\begin{bmatrix} x & y \end{bmatrix}
(V_1+V_1^T)
\begin{bmatrix} x \\ y\end{bmatrix}
\\ &= 2pq\cdot x^2+2(1-2kp)\cdot xy+2n\cdot y^2 
\\ &= \frac1{2pq}\left( \left(2pq\cdot x+(1-2kp)\cdot y\right)^2
+\left(4pqn-(2kp-1)^2\right)y^2  \right),\\
Q_0(x,y) &=
\begin{bmatrix} x & y \end{bmatrix}
(V_0+V_0^T)
\begin{bmatrix} x \\ y\end{bmatrix}
\\ &= 2pq\cdot x^2+2(1-2kp+2qr)\cdot xy+2(rs-2kr+n)\cdot y^2 
\\ &= \frac1{2pq}\left( \left(2pq\cdot x+(1-2kp+2qr)\cdot y\right)^2
+\left(4pqn-(2kp-1)^2\right)y^2  \right),
\end{align*}
where we used that 
\begin{gather*}
(1-2kp+2qr)^2-(2kp-1)^2 = 2qr(2-4kp+2qr)= \\
= 2qr(2ps-4kp) = 4pq(rs-2kr).
\end{gather*}

By \cref{thm:uv0v1}, the forms $2pqQ_0$ and $2pqQ_1$ are not isomorphic if
\begin{gather*}
1-2kp \nc 1-2kp+2qr \pmod {2pq},\\
1-2kp \nc -(1-2kp+2qr) \pmod {2pq},\\
4pqn-(2kp-1)^2 > \frac{(pq-1)^2}{3}.
\end{gather*}
The first condition is equivalent to $0 \nc r \pmod {p}$, which holds because $p>1$ and $ps-qr=1$.
The second condition is $0 \nc -2+4kp-2qr \pmod {2pq}$, which holds by our assumption $2kp\nc 1\pmod q$ (in fact, equivalent to it).
The third condition is equivalent to
\begin{gather*}
4pqn-(2kp-1)^2 \ge \frac{(pq-1)^2+2}{3}, \\
n \ge \frac{(2kp-1)^2}{4pq} + \frac{(pq-1)^2+2}{12pq} = \frac{k(pk-1)}q + \left( \frac{pq}{12} - \frac16 + \frac{1}{2pq}\right),
\end{gather*}
which was also an assumption. Therefore, the $Q_i$ are not isomorphic over $\Z$.

Moreover, $Q_0$ and $-Q_1$ are not isomorphic either, as the $Q_i$ are positive definite. Hence, the $X_i$ are not homeomorphic, thus the pairs $(D^4,\S_i)$ are not homeomorphic, and $\S_0\niso \S_1$.
\end{proof}

\begin{proof}[Proof of \cref{thm:pq}]
If $|p|=1$ or $|q|=1$, then the two annuli are easily seen to be isotopic in $S^3$ (even smoothly isotopic).

If $|p|,|q|\ne1$, then we can assume $p,q>1$. Suppose $A_0(p,q)\iso A_1(p,q)$. Let $k=0$, then $q>1$ implies $2kp\nc 1\pmod q$, and let $n$ be large enough (in fact, $n=1$ works). Then attaching the band $B(k,n)$ to $A_i(p,q)$, we get the surfaces $\S_i(p,q,k,n)=\S_i$. Viewing $D^4$ as $(S^3\times[0,1])\cup D^4_{1/2}$, we can isotope $\S_i$ into
\[
(B(k,n)\times\{0\})\cup (\d A_i(p,q)\times[0,1])\cup (A_i(p,q)\times\{1\}),
\]
and then do the isotopy $A_0(p,q)\iso A_1(p,q)$ in $D^4_{1/2}$.
We get $\S_0\iso \S_1$, contradicting \cref{thm:pqkn}. Hence $A_0(p,q)\niso A_1(p,q)$.
\end{proof}

\begin{proof}[Proof of \cref{thm:isotopic}]
Conway--Powell proved in \cite{Conway-Powell}*{Theorem 1.9} that for an Alexander polynomial one knot, any two Seifert surfaces of the same genus are topologically isotopic in the 4-ball. In our case, using the Seifert matrix of $\S_1$ from above, its boundary $K$ has Alexander polynomial
\begin{gather*}
\Delta_K(t) = \det(tV_1-V_1^T) =
\begin{vmatrix}
tpq-pq & -tkp-(1-kp) \\
t(1-kp)+kp & tn-n
\end{vmatrix} = \\
=p(nq-k(pk-1))(t-1)^2+t.
\end{gather*}
Hence, for $n=\frac{k(pk-1)}q$, we have $\Delta_K(t)=t\doteq 1$. Therefore, $\S_0\iso\S_1$.
\end{proof}

\begin{proof}[Proof of \cref{thm:n=0}]
When $n=0$, the quadratic forms can be factorized:
\begin{align*}
Q_0(x,y) &= 2\cdot \left(px+ry\right)\left(qx+\left(s-2k\right)y\right), \\
Q_1(x,y) &= 2\cdot x\left(pqx+\left(1-2kp\right)y\right).
\end{align*}

If $Q_0$ and $Q_1$ are isomorphic over $\Z$, then $Q_0(ax+by,cx+dy)=Q_1(x,y)$ for some invertible integral matrix $P=\begin{bmatrix} a & b \\ c & d \end{bmatrix}$.
Two factorized polynomials agree if and only if their linear factors agree up to multiplication by a constant, hence
\[
p(ax+by)+r(cx+dy)=\lambda x\quad \text{or}\quad \lambda\left(pqx+\left(1-2kp\right)y\right)
\]
for some real number $\lambda\ne 0$. Writing this in matrix form,
\[
\begin{bmatrix} a & c \\ b & d \end{bmatrix}
\begin{bmatrix} p \\ r \end{bmatrix}
=\lambda \begin{bmatrix} 1 \\ 0 \end{bmatrix}
\quad \text{or}\quad
\lambda \begin{bmatrix} pq \\ 1-2kp \end{bmatrix}.
\]
Combining this with the similar equation for the other linear factor, we get
\[
P^TM_0=M_1L_j\quad \text{for $j=0$ or 1},
\]
where
\begin{gather*}
P^T=\begin{bmatrix} a & c \\ b & d \end{bmatrix}, \quad
M_0=\begin{bmatrix} p & q \\ r & s-2k \end{bmatrix}, \quad
M_1=\begin{bmatrix} 1 & pq \\ 0 & 1-2kp \end{bmatrix}, \\
L_0=\begin{bmatrix} \lambda & 0 \\ 0 & 1/\lambda \end{bmatrix},
\quad \text{and}\quad
L_1=\begin{bmatrix} 0 & 1/\lambda \\ \lambda & 0 \end{bmatrix}.
\end{gather*}

If $Q_0$ and $-Q_1$ are isomorphic, then similarly, we get that $P^TM_0=M_1L_j$ for $j=2$ or 3, where
\[
L_2=\begin{bmatrix} \lambda & 0 \\ 0 & -1/\lambda \end{bmatrix},
\quad \text{and}\quad
L_3=\begin{bmatrix} 0 & -1/\lambda \\ \lambda & 0 \end{bmatrix}.
\]

Hence, we want to show that for large enough $k$ with $2kp\nc 1 \pmod q$, there do not exist a $j$ and a real number $\lambda\ne 0$ such that $P^T=M_1L_jM_0^{-1}$ is integral and invertible over $\Z$. Assume that there exist such a $j$ and $\lambda>0$ (without loss of generality).

If $j=0$, then
\[
P^TM_0=M_1L_0=\begin{bmatrix}
\lambda & pq/\lambda \\
0       & (1-2kp)/\lambda
\end{bmatrix}
\]
is integral, hence $\lambda$ is an integer. Since
\[
(P^T)^{-1}M_1=M_0L_0^{-1}=\begin{bmatrix}
p/\lambda & q \lambda \\
r/\lambda & (s-2k) \lambda
\end{bmatrix}
\]
is integral, we get $\lambda\mid p$ and $\lambda\mid r$, but $(p,r)=1$ then implies $\lambda=1$. Thus
\[
P^T=M_1L_0M_0^{-1}=\begin{bmatrix}
\dfrac{2k-s+pqr}{2kp-1} & \dfrac{q(1-p^2)}{2kp-1} \\
-r                      & p
\end{bmatrix}
\]
can only be integral if $2kp-1\mid q(1-p^2)$. Since $p>1$, this does not hold for $|k|$ large enough.

If $j=1$, then 
\[
P^TM_0=M_1L_1=\begin{bmatrix}
pq\lambda      & 1/\lambda \\
(1-2kp)\lambda & 0
\end{bmatrix}
\]
is integral, hence $\mu=1/\lambda$ is an integer. Then
\[
P^T=M_1L_1M_0^{-1}=\begin{bmatrix}
\dfrac{2kpq-pqs+r\mu^2}{(2kp-1)\mu} & \dfrac{pq^2-p\mu^2}{(2kp-1)\mu} \\
\dfrac{s-2k}{\mu}                   & \dfrac{-q}{\mu}
\end{bmatrix}
\]
is integral.
The bottom right entry implies $\mu\mid q$. The top right entry implies $2kp-1\mid pq^2/\mu-p\mu$. Since $(2kp-1,p)=1$, we have 
\[
2kp-1\mid q^2/\mu - \mu.
\]
If $|k|$ is large enough, then $|2kp-1| \ge q^2$. Since $|q^2/\mu - \mu|\le q^2-1$, the only possibility is $q^2/\mu - \mu = 0$, hence $\mu=q$. But the bottom left entry in $P^TM_0$ implies $\mu\mid 1-2kp$, contradicting $2kp\nc 1 \pmod q$.

If $j=2$, then the same argument as for $j=0$ shows that we need to have $2kp-1\mid q(1+p^2)$, which does not hold for $|k|$ large enough.

If $j=3$, then the same argument as for $j=1$ shows that $\mu=1/\lambda$ is an integer, $\mu\mid q$ and $2kp-1\mid q^2/\mu+\mu$, which does not hold if $|2kp-1|\ge q^2+2$.

Hence $Q_0$ is not isomorphic to $\pm Q_1$ over $\Z$ if $2kp\nc 1 \pmod q$ and
\[
|k|>k_0=\max\left(\frac{q(1+p^2)+1}{2p}, \frac{q^2+2}{2p}\right).
\]
Then the pairs $(D^4,\S_0)$ and $(D^4,\S_1)$ are not homeomorphic, and $\S_0\niso\S_1$.
\end{proof}

\begin{proof}[Proof of \cref{thm:cables}]
The proof is very similar to \cref{thm:pq}, so we only highlight the differences.

$A_0\cup A_1$ is a knotted torus $T$ in $S^3$ that bounds a solid torus. We attach a band $B(k,n)$ to $A_0$ and $A_1$ that is disjoint from the inside of the solid torus and has linking number $k$ with the core of the solid torus and has self-linking number $n$. We call the resulting surfaces $\S_0$ and $\S_1$, and will show that their Seifert matrices are the same as in the above proof.

The $(2p,2q)$-cable of the knot $K$ is defined using the longitude $\lambda$ of $T$ that satisfies $\lk(\lambda,\lambda^+)=0$ and a meridian $\mu$. Then the core $\alpha_1$ of $A_1$ is $p\lambda+q\mu$ in $H_1(T)$, hence
\[
\lk(\alpha_1, \alpha_1^+) = \lk(p\lambda, p\lambda^+ + q\mu^+) = pq,
\]
where $\lambda^+, \mu^+$ denote the push-off away from the solid torus. A base for $H_1(\S_i)$ is $\alpha_i,\beta_i$ similarly as in \cref{fig:pqkn2}. Then the Seifert matrix of $\S_1$ is the same $V_1$ as above. For the Seifert matrix of $\S_0$, we now have that $\beta_0-\beta_1=t_{r,s}$ is an $(r,s)$-cable of $K$ with $ps-qr=1$. Hence we have
\begin{align*}
\lk(t_{r,s},t_{r,s}^+) &= rs, \\
\lk(t_{r,s}, \alpha_0^+) &= \lk(r\lambda^+ + s\mu^+, p\lambda) = ps,
\end{align*}
and the rest of the calculation for $V_0$ is the same.

Therefore, for $p,q>1$ and suitably chosen $k,n$, the symmetrized Seifert forms are not isomorphic over $\Z$, hence $\S_0\niso\S_1$, and $A_0\niso A_1$.
\end{proof}

We end this section with a question. Note that \cref{thm:pq} was stated as an equivalence, but \cref{thm:cables} was stated as an implication. If $|p|=1$, then clearly $A_0\iso A_1$ for any $(2p,2q)$-cable. Is it true that if $|q|=1$, then $A_0\iso A_1$?

\section{Conway's algorithm}\label{sec:2}

There is a diagrammatic way to think about quadratic forms in two variables due to Conway \cite{Conway}. This can be used to solve generalized Pell's equations simply by looking at the corresponding diagram and it also gives an algorithm that decides whether two quadratic forms in two variables are isomorphic over $\Z$. We describe this algorithm and we will use it to study the symmetrized Seifert forms of the surfaces $\S_i(p,q,k,n)$ with values of $p,q,k,n$ that are not covered by \cref{thm:pqkn,thm:isotopic,thm:n=0}.

\subsection{Primitive vectors}
We begin the theory by discussing the structure of primitive vectors of $\Z^2$ and how it is related to an infinite 3-regular tree.

We call the vector $\v=(x,y)$ \emph{primitive} if $x,y$ are coprime. We say $(\v_1,\v_2)$ is a \emph{base} if their integral linear combinations give all of $\Z^2$ (then $\v_i$ must be primitive). We call $(\v_1,\v_2,\v_3)$ a \emph{superbase} if $\v_1+\v_2+\v_3=0$ and $(\v_1,\v_2)$ is a base. Then $(\v_1,\v_3)$ and $(\v_2,\v_3)$ is a base as well.

For a vector $\v$, we call the set $\pm \v=\{\v,-\v\}$ a \emph{lax vector}. If $(\v_1,\v_2)$ is a base, then we say $\{\pm\v_1,\pm\v_2\}$ is a \emph{lax base}, and if $(\v_1,\v_2,\v_3)$ is a superbase, then we say $\{\pm\v_1,\pm\v_2,\pm\v_3\}$ is a \emph{lax superbase}.

Note that each lax superbase $\{\pm\v_1,\pm\v_2,\pm\v_3\}$ contains 3 lax bases, and each lax base $\{\pm\v_1,\pm\v_2\}$ is contained in 2 lax superbases, given by $\pm\v_3=\pm(\v_1+\v_2)$ and $\pm\v_3=\pm(\v_1-\v_2)$. Hence the lax bases and lax superbases of $\Z^2$ can be visualized as a 3-regular graph, where the vertices are the lax superbases, and the edges are the lax bases, with an edge connecting to a vertex if the lax base is contained in the lax superbase. Conway then proved \cite{Conway} that this graph is connected and contains no cycles, hence it is an infinite 3-regular tree.

\begin{figure}
\par{\centering\def\svgwidth{0.85\textwidth}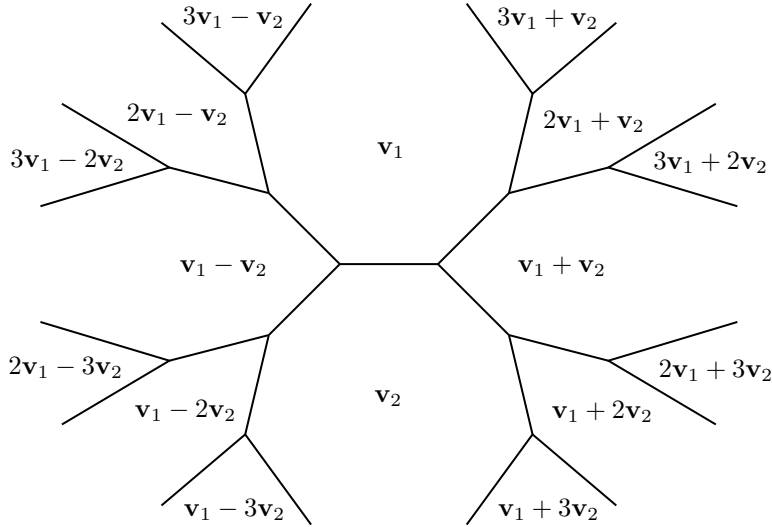\par}
\caption{The primitive vectors as the regions of a 3-regular tree, for $\v_1,\v_2$ a base of $\Z^2$.}\label{fig:conway1}
\end{figure}

If $\pm\v_1$ is a primitive lax vector, then it is in a lax base $\{\pm\v_1,\pm\v_2\}$, which is contained in a lax superbase $\{\pm\v_1,\pm\v_2,\pm\v_3\}$. From any such vertex that contains $\pm\v_1$, there are exactly 2 edges that contain $\pm\v_1$, namely $\{\pm\v_1,\pm\v_2\}$ and $\{\pm\v_1,\pm\v_3\}$. Hence we get an infinite path of edges that contain $\pm\v_1$. By drawing the tree in a suitable way on the plane, we can regard primitive lax vectors as `regions' of the tree that are bounded by these infinite paths. See \cref{fig:conway1}, where the $\pm$ signs have been omitted.

\subsection{The values of a quadratic form}
Consider an integral quadratic form in two variables: $Q(x,y)=Ax^2+Hxy+By^2$ with $A,B,H\in\Z$. We wish to understand the structure of values of $Q$ on $\Z^2$.

Since $Q(\lambda\v)=\lambda^2 Q(\v)$ for any integer $\lambda$, and $Q(-\v)=Q(\v)$, it is enough to understand the values of $Q$ on primitive lax vectors. We will write the value $Q(\v)$ on the region corresponding to $\pm\v$. We can start by writing $A$, $B$ and $A+H+B$ on the regions corresponding to $(1,0)$, $(0,1)$ and $(1,1)$; these regions have a common vertex. Then we can continue filling in the diagram by using the following observation: direct calculation shows that for any two vectors $\v_1,\v_2$,
\[
Q(\v_1+\v_2)+Q(\v_1-\v_2)=2(Q(\v_1)+Q(\v_2)).
\]
This means that if the 4 values around an edge are $a=Q(\v_1)$, $b=Q(\v_2)$, $c=Q(\v_1+\v_2)$ and $d=Q(\v_1-\v_2)$, then $d$, $a+b$, $c$ form an arithmetic progression.

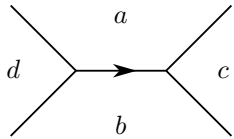
\begin{figure}[ht]
\par{\centering\def\svgwidth{0.25\textwidth}
\begingroup%
  \makeatletter%
  \providecommand\color[2][]{%
    \errmessage{(Inkscape) Color is used for the text in Inkscape, but the package 'color.sty' is not loaded}%
    \renewcommand\color[2][]{}%
  }%
  \providecommand\transparent[1]{%
    \errmessage{(Inkscape) Transparency is used (non-zero) for the text in Inkscape, but the package 'transparent.sty' is not loaded}%
    \renewcommand\transparent[1]{}%
  }%
  \providecommand\rotatebox[2]{#2}%
  \newcommand*\fsize{\dimexpr\f@size pt\relax}%
  \newcommand*\lineheight[1]{\fontsize{\fsize}{#1\fsize}\selectfont}%
  \ifx\svgwidth\undefined%
    \setlength{\unitlength}{183.23925829bp}%
    \ifx\svgscale\undefined%
      \relax%
    \else%
      \setlength{\unitlength}{\unitlength * \real{\svgscale}}%
    \fi%
  \else%
    \setlength{\unitlength}{\svgwidth}%
  \fi%
  \global\let\svgwidth\undefined%
  \global\let\svgscale\undefined%
  \makeatother%
  \begin{picture}(1,0.55406455)%
    \lineheight{1}%
    \setlength\tabcolsep{0pt}%
    \put(0,0){\includegraphics[width=\unitlength,page=1]{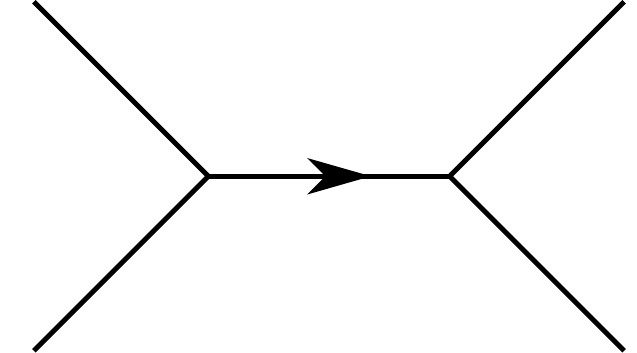}}%
    \put(0.51647819,0.47733667){\makebox(0,0)[t]{\lineheight{1.25}\smash{\begin{tabular}[t]{c}$a$\end{tabular}}}}%
    \put(0.51647819,0.01154359){\makebox(0,0)[t]{\lineheight{1.25}\smash{\begin{tabular}[t]{c}$b$\end{tabular}}}}%
    \put(0.94584767,0.24535645){\makebox(0,0)[t]{\lineheight{1.25}\smash{\begin{tabular}[t]{c}$c$\end{tabular}}}}%
    \put(0.06526427,0.24535645){\makebox(0,0)[t]{\lineheight{1.25}\smash{\begin{tabular}[t]{c}$d$\end{tabular}}}}%
  \end{picture}%
\endgroup%
\par}
\caption{$d<a+b<c$ is an arithmetic progression.}\label{fig:conway2}
\end{figure}

We also give an orientation to edges denoting the direction of increase of the arithmetic progression (if it is non-constant). Using the arithmetic progression rule, we can complete the labelling of all the regions and orientations of the edges, resulting in the diagram $\T(Q)$, called the \emph{topograph} of $Q$.

Two quadratic forms $Q_0$ and $Q_1$ are called isomorphic over $\Z$ if there exists an invertible integral matrix $P$ such that $Q_0(\v)=Q_1(P\v)$ for all $\v$. In other words, they represent the same function on $\Z^2$ but with respect to different bases. Hence, $Q_0$ and $Q_1$ are isomorphic over $\Z$ if and only if there is a vertex in $\T(Q_0)$ whose surrounding 3 regions have the same numbers as a vertex in $\T(Q_1)$. This shows that the topographs can be used to decide whether two forms are isomorphic over $\Z$, but in order to turn this into an algorithm, we will need to examine the structures of the topographs more closely.

\subsection{Classifying forms by signs}
\begin{figure}
\par{\centering\def\svgwidth{1\textwidth}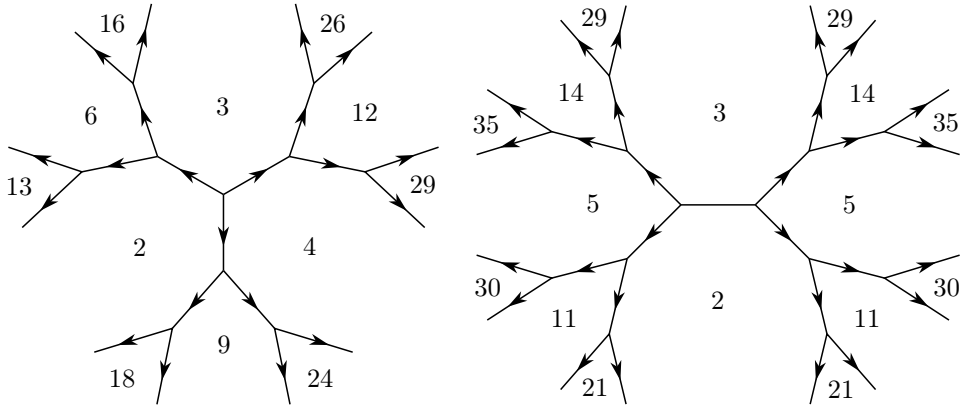\par}
\caption{Topograph of $2x^2+xy+3y^2$ (simple well) and $2x^2+3y^2$ (double well).}\label{fig:well}
\end{figure}

\begin{figure}
\par{\centering\def\svgwidth{1\textwidth}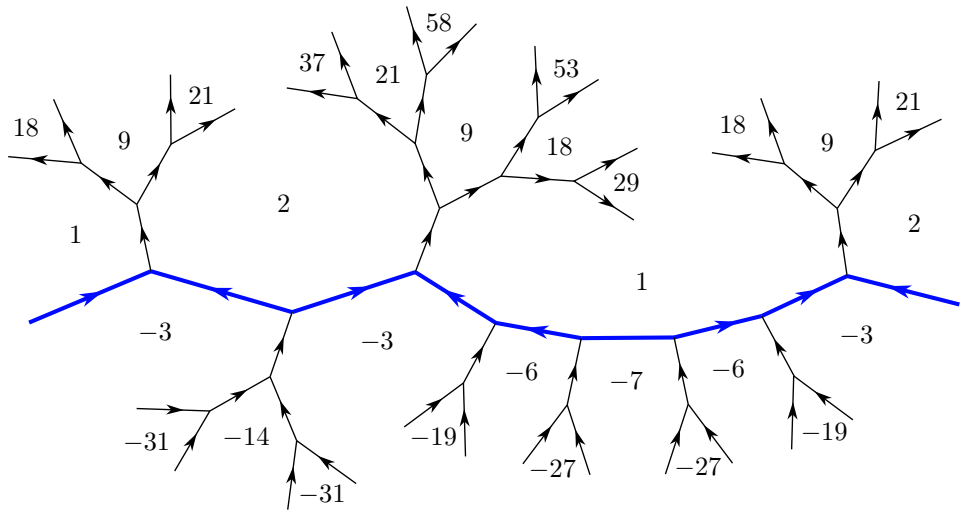\par}
\caption{Topograph of $x^2-7y^2$ (river).}\label{fig:river}
\end{figure}

\begin{figure}
\par{\centering\def\svgwidth{1\textwidth}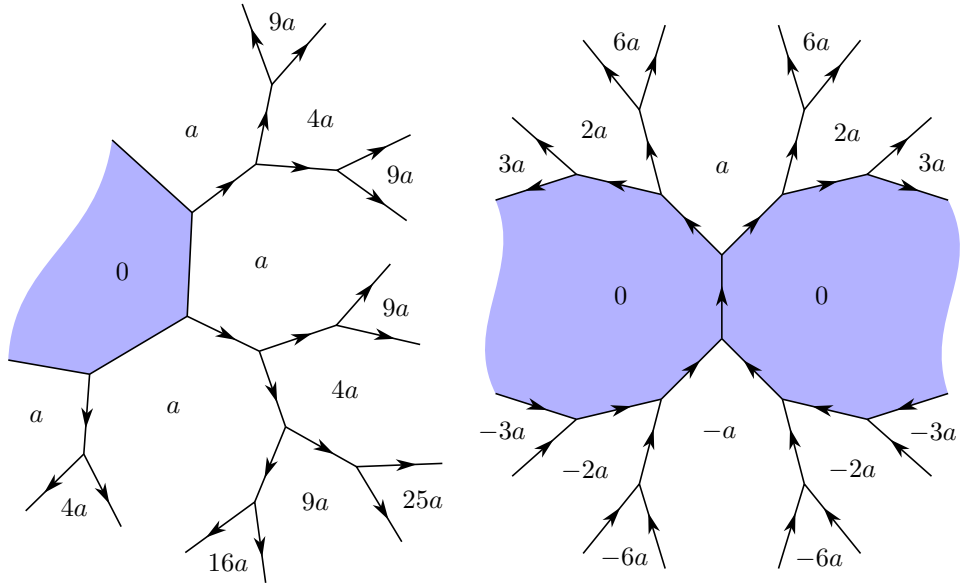\par}
\caption{Topograph of $ax^2$ (lake) and $axy$ (weir) for $a>0$.}\label{fig:lake-weir}
\end{figure}

\begin{figure}
\par{\centering\def\svgwidth{1\textwidth}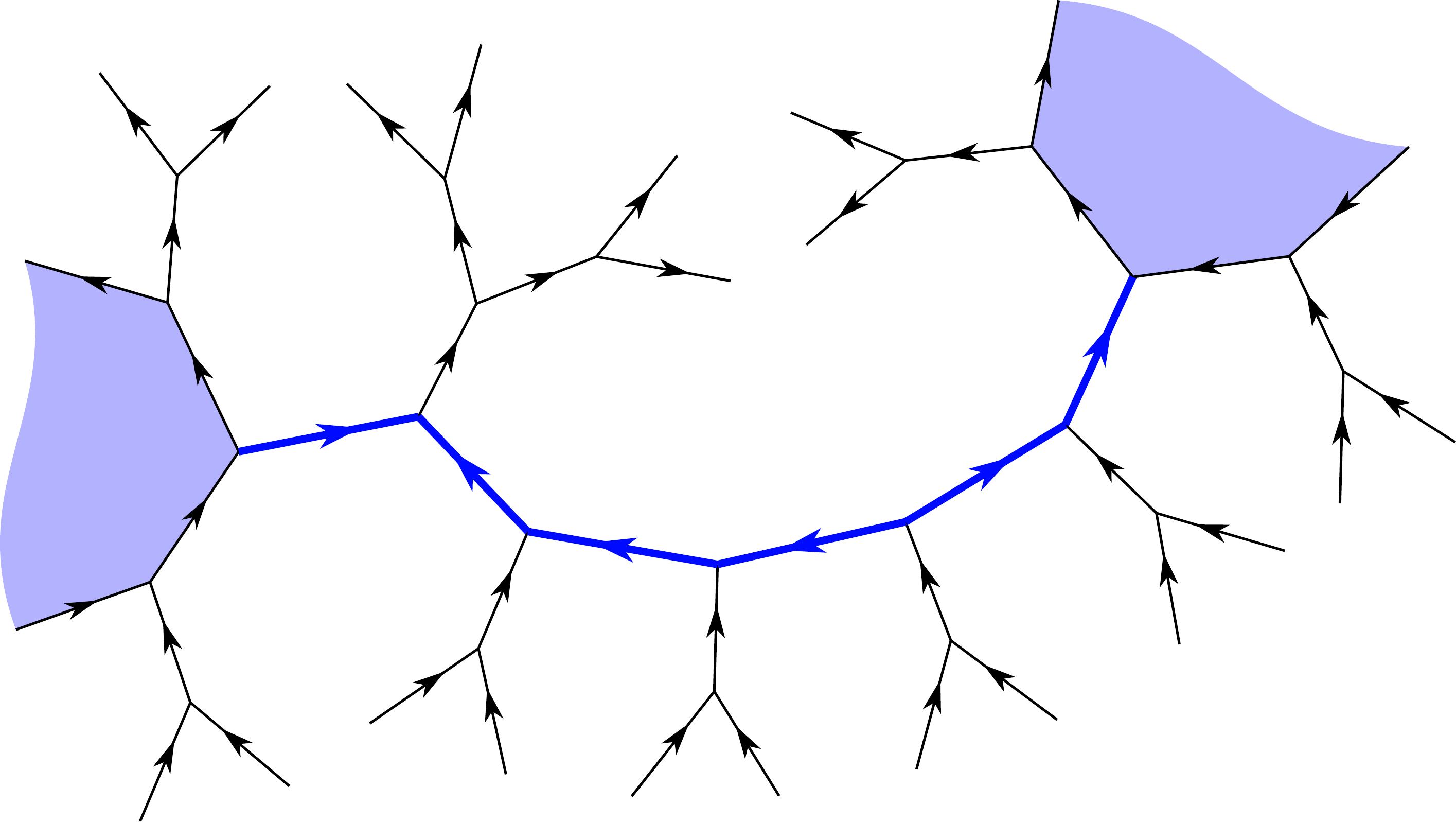\par}
\caption{Topograph of $6x^2+11xy$ (lake-pair).}\label{fig:doublelake}
\end{figure}

We get different types of structures depending on the signs of the values that a quadratic form takes. See \cref{fig:well,fig:river,fig:lake-weir,fig:doublelake} for examples. Conway proved the following:

\begin{itemize}
\item $+$ forms: If $Q$ is positive definite, then $\T(Q)$ has a vertex such that all edges in the tree are directed away from that vertex, or $\T(Q)$ has two adjacent vertices with a non-oriented edge between them such that all other edges are directed away from them. These 1 or 2 vertices are surrounded by the 3 or 4 smallest values in the topograph, hence we call this a (simple or double) \emph{well}.
\item $-$ forms: The negative definite forms are similar.
\item $+-$ forms: For indefinite forms, we consider the edges that separate a positive and a negative value. If 0 is not taken as a value, then every vertex has 0 or 2 such edges, hence we get an infinite path, which we call a \emph{river}. There is only one river, and the other edges are directed away from the river on its positive side and towards the river on its negative side. Moreover, for integral forms, the river is periodic, i.e.~there exists an $m$ such that moving along the river $m$ steps, we get the same values in the regions next to the river. This is the key fact that allows us to turn this into an algorithm.
\item $0$ form: For the 0 form, all the values in the topograph are 0.
\item $0+$ forms: For positive semi-definite forms (that are not the 0 form), there is exactly one region labelled 0, called the \emph{lake}, and all the regions adjacent to the lake have the same value. Moreover, all the edges not in the boundary of the lake are directed away from the lake.
\item $0-$ forms: The negative semi-definite forms are similar.
\item $0{+}{-}$ forms: For indefinite forms that take 0 as a value, there are two lakes, with a finite river connecting the two lakes. The positive and negative values in the topograph are separated by the river. All the edges not on the river are directed away from the river on its positive side and towards the river on its negative side. When the river has length 0, the lakes are adjacent, and we also call this a \emph{weir}.
\end{itemize}

\subsection{The algorithm}
Using the observations in the previous section, we can give an algorithm that decides whether two integral quadratic forms in two variables are isomorphic over $\Z$. In fact, we calculate a complete invariant of quadratic forms, i.e.~an invariant that is the same for two forms if and only if the forms are isomorphic over $\Z$.

We start with a triplet of numbers, the values of the form on a superbase, and write these numbers on the regions around a vertex in a 3-regular tree. If the signs of these three numbers are not all the same, then we have found a river or a lake. Otherwise, e.g.~if they are all positive, then we calculate the orientations of the three edges using the arithmetic progression rule. (If the numbers are all negative, then we multiply everything by $(-1)$ and refer to the positive case.) If all edges are directed away from the vertex, then we have found a well. Otherwise, there is an edge directed towards the vertex, and we move along that edge and repeat these steps with the triplet around the new vertex. Continuing this procedure of descending, eventually we will find a well, a river, or a lake. If we have a river, then walking along the river, we either get to a lake or find that the river is periodic. If we have a lake, then descending along the boundary of the lake, we either get to a river or find that all the values are the same next to the lake. In each case, using the values around a well, river, or lake, we can clearly define a complete invariant.

We give an implementation of this algorithm in Wolfram Mathematica 13.3. In the following code, \texttt{invariant[tri\_]} calculates a complete invariant of integral quadratic forms in two variables. It takes a triplet of numbers $\mathtt{tri}=\{a,b,c\}$ as input (the values of the form on a superbase), and it gives an output that has the format
\{"WELL", $\{a, b, c\}$\},
\{"RIVER", $\{a, b, c\}$\},
\{"LAKE", $\{a\}$\},
\{"LAKE-PAIR", $\{\{a, b\}, \{c, d\}\}$\}, or
\{"WEIR", $\{a\}$\},
with some integers $a,b,c,d$ that are canonical in some sense (e.g.~smallest triplet along the river).

\pagebreak
\begin{scriptsize}
\begin{verbatim}
lastpos[list_] := Block[{i = -2}, While[list[[i]] < 0, i = i - 1]; list[[i]]]
lastneg[list_] := Block[{i = -2}, While[list[[i]] > 0, i = i - 1]; list[[i]]]
triplet[list_] := {Last@list, lastpos@list, lastneg@list}
rivernext[list_] := (
  {last, lastpositive, lastnegative} = triplet[list];
  Append[list, If[Positive[last],
    2 (lastnegative + last) - lastpositive, 
    2 (lastpositive + last) - lastnegative]])
river[tri_] := 
  NestWhile[rivernext, tri, Last@# != 0 &&
    (Length@# <= 4 || {last, lastpositive, lastnegative} != triplet[tri]) &]
riverinvariant[river_] := 
  Block[{t = Table[Sort@triplet[river[[;; i]]], {i, 3, Length@river - 2}]}, 
    t[[First@Ordering[t, 1]]]]
descendnext[tri_] := Block[{min, middle, max},
  {min, middle, max} = Sort[tri];
  If[max > min + middle,
    {min, middle, 2 (min + middle) - max},
    {"WELL", {min, middle, max}}]]
descend[tri_] := NestWhile[descendnext, tri, Length@# == 3 && Min@# > 0 &]
invariant[tri_] := Block[{d = tri, r, min, middle, max, a, b, u, v},
  If[Min[tri] > 0, d = descend[tri]; If[First@d == "WELL", Return[d]]];
  If[Max[tri] < 0, d = -descend[-tri]; 
    If[First@d == -"WELL", Return[{"WELL", Last@d}]]];
  If[! MemberQ[d, 0],
    {min, middle, max} = Sort[d];
    r = river[{min, max, middle}];
    If[Last@r == 0, d = triplet[r], Return[{"RIVER", riverinvariant[r]}]]];
  {a, b} = Sort@Rest@SortBy[d, Abs];
  If[a == b, Return[{"LAKE", {a}}]];
  {a, b} = {Mod[a, b - a, -b + a + 1], Mod[a, b - a, 1]};
  If[a == 0, Return[{"WEIR", {b}}]];
  r = river[{a, b, 2 (a + b)}];
  {u, v} = Sort@DeleteCases[triplet[r], 0];
  Return[{"LAKE-PAIR", Sort@{{a, b}, {u, v}}}]]
\end{verbatim}
\end{scriptsize}

\subsection{Application to the Seifert forms}\label{sec:discussion}

We can use the algorithm in the previous section on our quadratic forms $Q_0, Q_1$ that arose as the symmetrized Seifert forms of the surfaces $\S_i(p,q,k,n)=\S_i$. For instance, Lyon's family of examples \cite{Lyon-example} is $\S_i(3,-4,1,-n)$ in our notation, which is the reflection of $\S_i(3,4,1,n)$. For $n=1$, this is not covered by \cref{thm:pqkn}, but Conway's algorithm can be used to show that the symmetrized Seifert forms are non-isomorphic over $\Z$, and thus Lyon's surface-pairs are non-isotopic in $D^4$ for $n>0$.

We can assume $p,q>0$. Recall that the matrix for $Q_i$ is $V_i+V_i^T$, where
\[
V_0=
\begin{bmatrix}
pq    & qr-kp\\
ps-kp & rs-2kr+n
\end{bmatrix},
\quad
V_1=
\begin{bmatrix}
pq   & -kp\\
1-kp & n
\end{bmatrix},
\quad
ps-qr=1.
\]
The following code calculates whether the invariants of $Q_0,Q_1$ agree for concrete values of $p,q,k,n$ (in the example below, $(p,q)=(2,3)$ is fixed and $|k|,|n|\le30$ varies), and plots a diagram of the result:

\begin{scriptsize}
\begin{verbatim}
{p, q, r, s, k, n} =.;
tri0 = {p q, p q + (p s + q r - 2 k p) + (r s - 2 k r + n), r s - 2 k r + n};
tri1 = {p q, p q + (1 - 2 k p) + n, n};
transform[k_, n_] := 
  Sort[{{Mod[k, q], n - (k (p k - 1))/q + (Mod[k, q] (p Mod[k, q] - 1))/q},
    {Mod[s - k, q], n - (k (p k - 1))/q + (Mod[s - k, q] (p Mod[s - k, q] - 1))/q}}]
parallelgetsol[kmin_, kmax_, nmin_, nmax_] := AbsoluteTiming[
  unknown = Flatten[Table[If[! MemberQ[Keys@sol, transform[k, n]], {k, n}, Nothing],
    {n, nmax, nmin, -1}, {k, kmin, kmax}], 1];
  paralleltable = ParallelTable[{k, n} = value;
    {k, n, invariant[tri0] != invariant[tri1]},
    {value, unknown}, ProgressReporting -> False];
  (sol[transform[#[[1]], #[[2]]]] = #[[3]]) & /@ paralleltable;]
normalize[sol_] := Block[{normalsol = <||>},
  Do[If[MemberQ[Keys@sol, transform[k, n]], 
    normalsol[{k, n}] = sol[transform[k, n]]], {k, -size, size}, {n, -size, size}];
  normalsol]
plot[sol_] := Block[{normalsol = normalize[sol]},
  Column[{{p, q},
    MatrixPlot[Normal@SparseArray[
      Table[{-x[[2]] + size + 1, x[[1]] + size + 1} -> 
        If[normalsol[x], 1, -1], {x, Keys@normalsol}],
      {2 size + 1, 2 size + 1}], FrameTicks -> None, ImageSize -> 250],
    Iconize[normalsol]}, Center]]

{p, q} = {2, 3};
s = ModularInverse[p, q];
r = (p s - 1)/q;
sol = <||>;
size = 30;
Table[parallelgetsol[k, k, -size, size], {k, -size, size}]; // AbsoluteTiming // First
plot[sol]
\end{verbatim}
\end{scriptsize}

\begin{figure}[p]
\par{\centering\def\svgwidth{0.98\textwidth}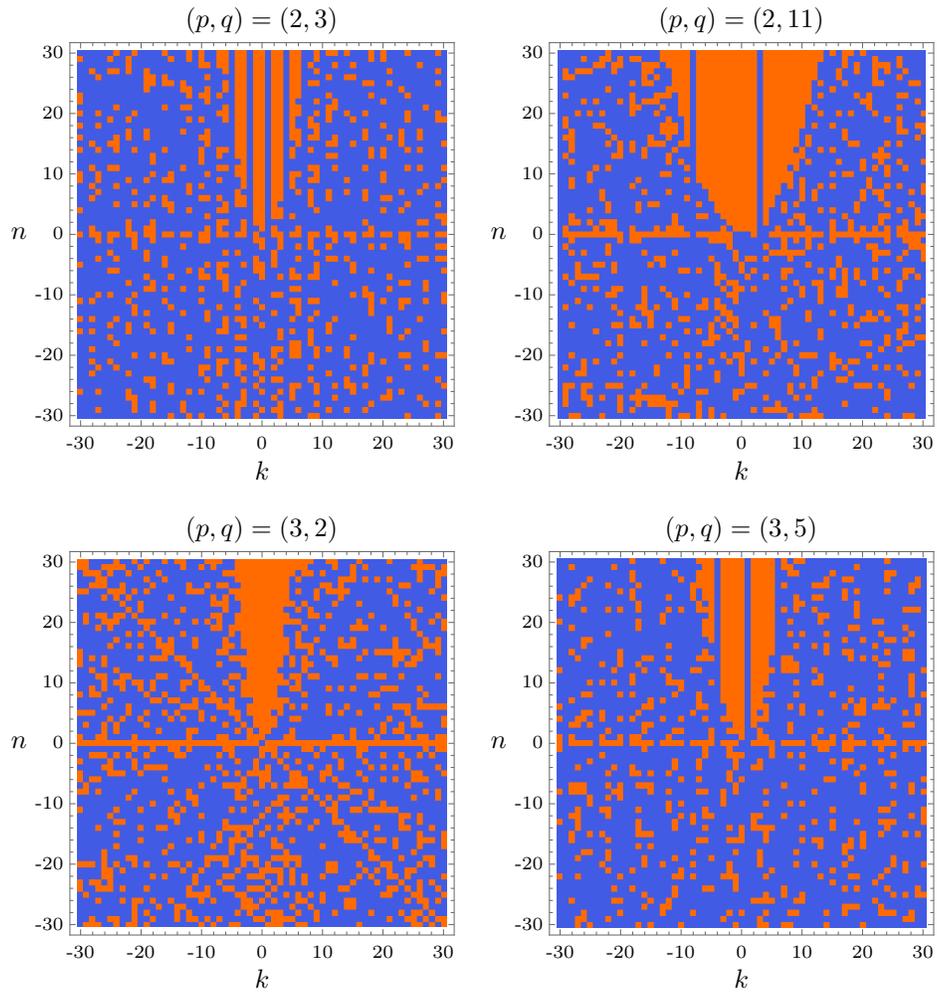\par}
\caption{Diagrams showing when the symmetrized Seifert forms are isomorphic (blue) and non-isomorphic (orange) over $\Z$.}\label{fig:parabolas}
\end{figure}

\begin{figure}[p]
\par{\centering\def\svgwidth{1\textwidth}
\begingroup%
  \makeatletter%
  \providecommand\color[2][]{%
    \errmessage{(Inkscape) Color is used for the text in Inkscape, but the package 'color.sty' is not loaded}%
    \renewcommand\color[2][]{}%
  }%
  \providecommand\transparent[1]{%
    \errmessage{(Inkscape) Transparency is used (non-zero) for the text in Inkscape, but the package 'transparent.sty' is not loaded}%
    \renewcommand\transparent[1]{}%
  }%
  \providecommand\rotatebox[2]{#2}%
  \newcommand*\fsize{\dimexpr\f@size pt\relax}%
  \newcommand*\lineheight[1]{\fontsize{\fsize}{#1\fsize}\selectfont}%
  \ifx\svgwidth\undefined%
    \setlength{\unitlength}{577.125bp}%
    \ifx\svgscale\undefined%
      \relax%
    \else%
      \setlength{\unitlength}{\unitlength * \real{\svgscale}}%
    \fi%
  \else%
    \setlength{\unitlength}{\svgwidth}%
  \fi%
  \global\let\svgwidth\undefined%
  \global\let\svgscale\undefined%
  \makeatother%
  \begin{picture}(1,0.36794673)%
    \lineheight{1}%
    \setlength\tabcolsep{0pt}%
    \put(0,0){\includegraphics[width=\unitlength,page=1]{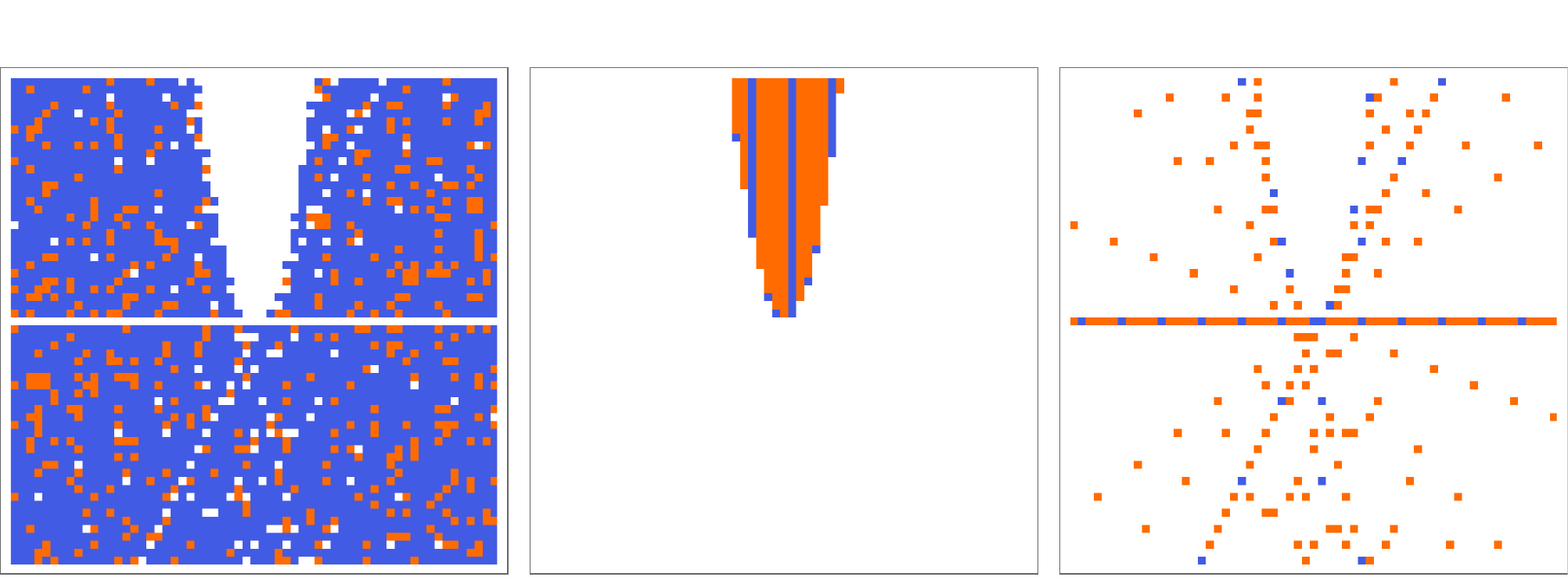}}%
    \put(0.16211826,0.33509261){\makebox(0,0)[t]{\lineheight{1.25}\smash{\begin{tabular}[t]{c}river ($+-$)\end{tabular}}}}%
    \put(0.5,0.33509261){\makebox(0,0)[t]{\lineheight{1.25}\smash{\begin{tabular}[t]{c}well ($+$)\end{tabular}}}}%
    \put(0.83788173,0.33509261){\makebox(0,0)[t]{\lineheight{1.25}\smash{\begin{tabular}[t]{c}lake-pair ($0{+}{-}$)\end{tabular}}}}%
  \end{picture}%
\endgroup%
\par}
\caption{Splitting the diagram by topograph type in the case $(p,q)=(3,5)$.}\label{fig:split}
\end{figure}
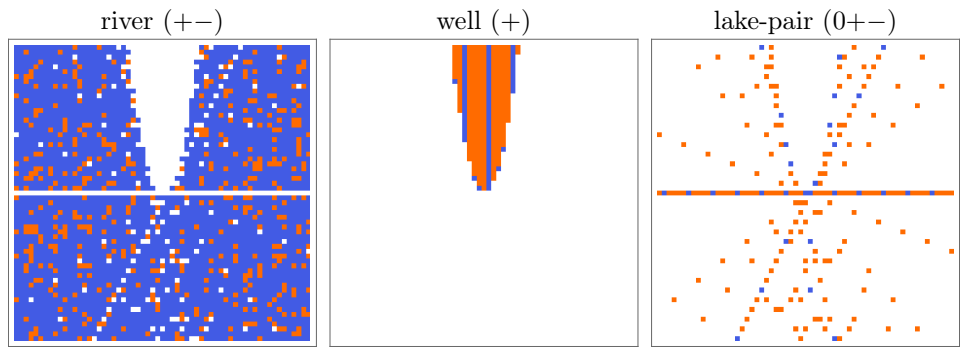

\cref{fig:parabolas} shows the resulting diagrams for a few different values of $(p,q)$. An orange square in these diagrams means that $Q_0$ and $Q_1$ are not isomorphic over $\Z$, hence the pairs $(D^4,\S_0)$ and $(D^4,\S_1)$ are not homeomorphic in an orientation-preserving way, and $\S_0\niso\S_1$. A blue square means that the Seifert form does not give an answer to whether the surfaces are isotopic. Apart from when the Alexander polynomial is 1 (the points on the parabola, see \cref{thm:isotopic}), we do not know if they are isotopic or not.

\subsubsection{Topograph type}
We can also split the diagrams by the type of the topographs to gain more insight, see \cref{fig:split}.

We show that the type of topograph (river, well, lake, or lake-pair) is the same for the two quadratic forms. $Q_0$ and $Q_1$ are always isomorphic over $\Q$, because
\[
P^T V_0P=V_1 \quad \text{for}\quad P=
\begin{bmatrix}
1 & -r/p \\
0 & 1
\end{bmatrix}.
\]
Hence the signs ($0{+}{-}$) of the values of $Q_0$ on $\Q^2$ are the same as of $Q_1$, thus their signs on $\Z^2$ agree as well.
It follows that $\T(Q_0)$ and $\T(Q_1)$ have the same type.

\subsubsection{Symmetries}
We note two important symmetries of these diagrams. Consider the affine transformations
\begin{align*}
\tau(k,n) &= (k+q, \; n+2kp+pq-1), \\
\rho(k,n) &= (s-k, \; n-2kr+rs).
\end{align*}
Note that the value of $n-\frac{k(kp-1)}{q}$ is invariant under $\tau$ and $\rho$. Intuitively, this means that the vertical distance from the parabola $n=\frac{k(kp-1)}{q}$ is unchanged, hence $\tau$ is translation to the right by $q$ `along the parabola', and $\rho$ is reflection in the vertical line $k=s/2$ along the parabola.

Direct calculation shows that the isomorphism classes of the quadratic forms $Q_0$ and $Q_1$ are invariant under $\tau$, while $\rho$ swaps them. Hence, the diagrams in \cref{fig:parabolas} are invariant under $\tau$ and $\rho$. In fact, the surface-pair $\S_i(p,q,k,n)$ is isotopic in $S^3$ to a surface-pair $\S_i(p,q,\tau(k,n))$ and to a $\S_{1-i}(p,q,\rho(k,n))$. This can be seen by sliding the attaching region of the band $B(k,n)$ along the full circle boundary of $A_i(p,q)$.

\subsubsection{Vertical lines}
The diagrams in \cref{fig:parabolas} show vertical blue lines $q$ distance apart, but only when $q$ is odd. Indeed, we prove that the lines $2kp\equiv 1 \pmod q$ are blue. The quadratic forms $(x,y)\mapsto(ux+v_iy)^2+ty^2$ for $i=0,1$ are clearly isomorphic over $\Z$ if $v_0\equiv -v_1 \pmod u$ (c.f.~\cref{thm:uv0v1}). We have seen that the forms $2pqQ_i$ are of this form with
\[
u=2pq, \quad v_0=1-2kp+2qr, \quad v_1=1-2kp, \quad t=4pqn-(2kp-1)^2.
\]
Hence, we need to show $v_0+v_1 = 2(1-2kp+qr) \equiv 0 \pmod {2pq}$. As $p,q$ are coprime, this holds if and only if $1-2kp+qr$ is divisible by $p$ and $q$. It is always divisible by $p$ because $1+qr=ps$, and it is divisible by $q$ when $2kp\equiv 1 \pmod q$ holds. Then $2pqQ_i$ are indeed isomorphic.

Another way to see this is by observing that every point $(k,n)$ on these vertical lines is fixed under a transformation $\tau^{2m}\circ\rho$ for some $m\in\Z$. Since $\tau$ fixes the isomorphism class of $Q_i$ and $\rho$ swaps them, we get
\[
Q_0^{(k,n)}=Q_0^{\tau^{2m}(\rho(k,n))}\cong Q_0^{\rho(k,n)}\cong Q_1^{(k,n)}.
\]

\subsubsection{The parabola}
The most prominent feature of the diagrams in \cref{fig:parabolas} is the parabolic region: we have proved in \cref{thm:pqkn} that $2kp\nc 1 \pmod q$ and
\[
n \ge \frac{k(pk-1)}q + \left( \frac{pq}{12} - \frac16 + \frac{1}{2pq} \right)
\]
gives a set of orange squares.
On the other hand, a consequence of \cref{thm:isotopic} is that the points along the parabola
$n = \frac{k(pk-1)}q$
are blue. Note that since $pq>0$, the quadratic forms $Q_i$ are positive definite if and only if $0<\det(V_1+V_1^T) = 4pqn-(2kp-1)^2$, or equivalently, $n>\frac{k(pk-1)}q$.

It is not true that $n > \frac{k(pk-1)}q$ and $2kp\nc 1 \pmod q$ gives a fully orange region, as the following example shows: let $(p,q,k,n)=(3,5,-1,1)$. Then
\[
V_0=\begin{bmatrix}
15 & 8 \\
9 & 5
\end{bmatrix} \quad \text{and} \quad
V_1=\begin{bmatrix}
15 & 3 \\
4 & 1
\end{bmatrix}
\]
yield isomorphic quadratic forms. The isomorphism $P^TV_0P=V_1$ is given by
$
P=\begin{bmatrix}
2 & 1 \\
-5 & -2
\end{bmatrix}.
$
In fact, the quadratic forms $2pqQ_1(x,y)$ and $2pqQ_0(x+y,-y)$ are the ones in \cref{rem:quad11} with $t=11$, which was the borderline case. See also the second diagram in \cref{fig:split}, which displays the positive definite forms, showing a few blue squares at the bottom of the region.

However, we prove that if $k=0$ or $q \le 3$, then the assumption on $n$ in \cref{thm:pqkn} can be weakened to $n > \frac{k(pk-1)}q$. For that, we use the improved bound on $t$ described in \cref{rem:t0}.

If $k=0$, then the forms $Q_0,Q_1$ have $u=2pq$, $v_0=1+2qr$, and $v_1=1$. Since $[v_0]_u\ge1=[v_1]_u$ and $[cv_0]_u^2\ge 0$ for all $c$, we have
\[
t_1 = \max_{c\ge 2\text{ integer}} \frac{[v_1]_u^2-[cv_0]_u^2}{c^2-1} \le
\max_{c\ge 2\text{ integer}} \frac{1}{c^2-1}=\frac13.
\]
We have seen that $[v_0]_u\ne[v_1]_u$ if $p>1$ and $2kp\nc1\pmod q$. Hence, $Q_0$ and $Q_1$ are non-isomorphic for every $0<t=4pqn-(2kp-1)^2$, which is what we wanted to prove.

If $q \le 3$, then we show that using the symmetries $\tau$ and $\rho$, it can be assumed that $k=0$. Indeed, if $q=2$, then $ps-qr=1$ implies that $s$ is odd, hence (using $\rho$) one of $k$ or $s-k$ is even, which can be translated (using $\tau$) to $0$ by some multiple of $q$. If $q=3$, then the condition $2kp\nc1\pmod q$ and $(q,s)=1$ implies $2k\equiv 2kps\nc s\pmod q$. Since $\rho$ is a reflection in the line $2k=s$, it means that it swaps the other two remainder classes modulo 3. Hence $k$ can be sent to $0$ by $\tau$ and $\rho$, and we are done.

\subsubsection{Horizontal and diagonal lines}
Lastly, the third diagram in \cref{fig:split} explains why the horizontal line $n=0$ was special in \cref{thm:n=0}. The diagram shows when the quadratic forms take 0 as a value, which means that each of the forms can be written as a product of two linear terms. This is what we used in the proof.

There often are other lines visible on this diagram too, but in the case of $(p,q)=(3,5)$, all the other lines are actually consequences of the line $n=0$ and the symmetries $\tau$ and $\rho$. For instance, the line $n=2k-2$ is the $\rho$-image of the line $n=0$ (with $(r,s)=(1,2)$). The line $n=-2k/5$ appears because the point $(5m,-2m)$ is the $\tau^{2m}$-image of the point $(-5m,0)$ on the $n=0$ line.

In other cases though, e.g.~for $(p,q)=(3,10)$, there are other factorizable quadratic forms that do not arise in this way.

\renewcommand{\eprint}[1]{\href{https://arxiv.org/pdf/#1}{arXiv:#1}}
\renewcommand{\MR}[1]{ \href{http://www.ams.org/mathscinet-getitem?mr=#1}{MR#1}}
\bibliography{bibliography}

\end{document}